\documentclass[12pt, a4paper]{article}

\usepackage[latin9]{inputenc}
\usepackage[T1]{fontenc}
\usepackage[french, english]{babel}
\usepackage{amsmath, amssymb, amsthm}
\usepackage{url}
\usepackage{hyperref}
\usepackage{geometry}
\usepackage{fancyhdr}
\newcounter{moncompteur}

\newtheorem{theorem}[moncompteur]{Théorème}
\newtheorem*{theoremnn}{Théorème}
\newtheorem{prop}[moncompteur]{Proposition}

\newtheorem{lemme}[moncompteur]{Lemme}

\renewcommand\o{\H{o}}
\newcommand\R{\mathbb{R}}
\newcommand\N{\mathbb{N}}
\newcommand\C{\mathbb{C}}
\newcommand\Z{\mathbb{Z}}
\renewcommand\O[1]{O\left( #1\right)}
\newcommand\MR{Matomäki et Radziwi\l\l\ }
\newcommand\MRen{Matomäki and Radziwi\l\l\ }
\newcommand\MRf{Matomäki et Radziwi\l\l}
\renewcommand\leq{\leqslant}
\renewcommand\geq{\geqslant}
\renewcommand\theta{\vartheta}
\newcommand\im{\Im\mathfrak{m}}
\newcommand\re{\Re\mathfrak{e}}
\newcommand\e{\mathrm{e}}
\renewcommand\epsilon{\varepsilon}

\begin{document}
\selectlanguage{french}
\title{Théorème d'Erd\o s-Kac\\dans presque tous les petits intervalles}
\author{Élie \bsc{GOUDOUT}}
\date{}
\maketitle
\selectlanguage{english}
\begin{abstract} 
We show that the Erd\o s-Kac theorem is valid in almost all intervals $\left[x,x+h\right]$ as soon as $h$ tends to infinity with $x$.  We also show that for all $k$ near $\log\log x$, almost all interval $\left[x,x+\exp\left(\left(\log\log x\right)^{1/2+\varepsilon}\right)\right]$ contains the expected number of integers $n$ such that $\omega(n)=k$. These results are a consequence of the methods introduced by \MRen to estimate sums of multiplicative functions over short intervals.
\selectlanguage{french} 
\begin{center} 
\textbf{Résumé} 
\end{center} 
\indent\indent On démontre que le théorème d'Erd\o s-Kac est valable dans presque tous les intervalles $[x,x+h]$ dès que $h$ tend vers l'infini avec $x$. On démontre aussi que pour tout $k$ proche de $\log\log x$, presque tout intervalle $\left[x,x+\exp\left(\left(\log\log x\right)^{1/2+\varepsilon}\right)\right]$ contient le nombre attendu d'entiers $n$ tels que $\omega(n)=k$. Ces résultats sont des conséquences des méthodes de \MR sur les fonctions multiplicatives dans les petits intervalles.
\end{abstract}
\selectlanguage{french}

\section{Introduction}

\indent\indent Pour $n\geq 1$, on note $\omega(n)$ le nombre de facteurs premiers distincts de $n$. Cette fonction admet la fonction $n\mapsto\log_2n$ comme ordre normal\footnote{Pour $k\geq 2$, $\log_k$ désigne la fonction $\log$ itérée $k$ fois. Par exemple, $\log_2 x=\log\log x$.}. En 1939, Erd\H{o}s et Kac montrent que la répartition de $\omega$ suit une loi normale. En 1958, Rényi et Tur\'an améliorent leur résultat de la manière suivante.
\begin{theoremnn}[Erd\H{o}s \& Kac~\cite{erdoskac} ; Rényi \& Tur\'an~\cite{renyituran}]
Uniformément pour $X\geq 20$ et $y\in\R$, on a
\[\frac{1}{X}\#\left\{ 1\leq n\leq X\,:\quad\omega(n)\leq\log_2 X+y\sqrt{\log_2 X}\right\}=\Phi(y)+\O{\frac{1}{\sqrt{\log_2 X}}},\]
où
\begin{equation}
\Phi(y):=\frac{1}{\sqrt{2\pi}}\int_{-\infty}^y\!\!\e^{-\tau^2/2}\mathrm{d}\tau.
\end{equation}
\end{theoremnn}
Comme précisé après le théorème \MakeUppercase{\romannumeral 3}.4.15 de~\cite{Tenenbaum}, le terme d'erreur, s'il est indépendant de $y$, est optimal. En 1981, Babu~\cite{babu} a montré que le même résultat (sans terme d'erreur) était valable dans tous les intervalles $\left(x,x+h\right]$ lorsque $y\in\R$ est fixé dès que $h\geq x^{a(x)/(\log_2 x)^{1/2}}$ où $a(x)\rightarrow +\infty$. On montre ici que, lorsque $h$ tend vers l'infini arbitrairement lentement avec $x$, le résultat se vérifie dans presque tous les intervalles $\left(x,x+h\right]$. Par ailleurs, Delange~\cite{delange} a montré en 1959 que le $\O{\left(\log_2 X\right)^{-1/2}}$ pouvait être remplacé par
\[\frac{\e^{-y^2/2}}{\sqrt{2\pi\log_2 X}}\left(\frac{2}{3}-c_1-\frac{y^2}{6}-\left\langle\log_2 X+y\sqrt{\log_2 X}\right\rangle\right)+\O{\frac{1}{\log_2 X}},\]
où
\begin{equation}\label{mertensconstant}
c_1:=\gamma+\sum_{p}\left(\log\left(1-\frac{1}{p}\right)+\frac{1}{p}\right)\simeq 0,\!261497
\end{equation}
est la constante de Mertens, et $\langle t\rangle$ désigne la partie fractionnaire de $t\in\R$. Ce résultat peut être retrouvé par les méthodes développées par Esseen dans sa thèse~\cite[p. 53-59]{theseEsseen}, avec un terme d'erreur légèrement plus faible. On s'en inspire donc pour retrouver l'analogue du résultat de Delange dans les petits intervalles. Pour $20\leq h\leq X$, $x\in\left[X,2X\right]$ et $y\in\R$, on définit alors
\begin{align}
F_{(x,x+h]}(y) &:=\frac{1}{h}\#\left\{x<n\leq x+h\,:\quad\omega(n)\leq\log_2 X +y\sqrt{\log_2 X}\right\},\\
\Phi_X(y) &:=\frac{1}{\sqrt{2\pi}}\int_{-\infty}^y\e^{-\tau^2/2}\mathrm{d}\tau+\frac{\e^{-y^2/2}}{\sqrt{2\pi\log_2 X}}\left(\frac{2}{3}-c_1-\frac{y^2}{6}-\left\langle\log_2 X+y\sqrt{\log_2 X}\right\rangle\right).\label{defPhi}
\end{align}
\begin{theorem}\label{petits intervalles}
Uniformément pour $20\leq h\leq X$ et $\alpha\geq 1$, on a
\[\left\Vert F_{(x,x+h]}-\Phi_X\right\Vert_{\infty}\ll\frac{\log_3 X}{\log_2 X}+\alpha^2\frac{\left(\log_2 h\right)^2}{\log h}\]
pour tout $x\in\left[X,2X\right]$ sauf sur un ensemble de mesure au plus
\[\ll X\left(\frac{1}{\left(\log h\right)^{\alpha}}+\frac{1}{\left(\log X\right)^{1/150}}\right).\]
\end{theorem}
\pagebreak[2]
Le paramètre $\alpha\geq 1$ pouvant varier uniformément, on peut diminuer la taille de l'ensemble exceptionnel au prix d'une légère perte sur le terme d'erreur, en prenant par exemple $\alpha=\log_2 h$. On note qu'il est possible, avec la méthode d'Esseen, d'augmenter la précision de $\Phi_X$ pour remplacer le $\O{\frac{\log_3 X}{\log _2 X}}$ par $o\left(\left(\log_2 X\right)^{-n/2}\right)$ pour tout $n\geq 1$ fixé. On renvoie le lecteur intéressé  à~\cite[p.~60-61]{theseEsseen}.\break Il semble par ailleurs possible d'obtenir $F_{(x,x+h]}(y)\sim\Phi(y)$ pour $y$ fixé, pour presque tout $x$, par la méthode des moments qui est plus élémentaire.

Pour $k\geq 1$, on note $\pi_k(x):=\#\left\{n\leq x\,:\quad\omega(n)=k\right\}$. La bonne qualité du terme d'erreur du Théorème~$\ref{petits intervalles}$ permet de retrouver $\pi_k(x+h)-\pi_k(x)$ pour des $k$ normaux et $h$ suffisamment grand. Il est cependant possible de faire cela en s'intéressant directement aux lois locales.
\begin{theorem}\label{loislocales}
Soient $\epsilon >0$ fixé et $X\geq 20$. Soient $h\leq X$ et $k$ un entier tels que
\begin{align}
\left\vert k-\log_2 X\right\vert\ll\sqrt{\log_2 X},\\
h\geq\exp\left((\log_2X)^{1/2+\epsilon}\right)\label{borneinfh}.
\end{align}
Alors, pour presque tout $x\sim X$, on a
\begin{equation}
\pi_k(x+h)-\pi_k(x)\sim \frac{h}{\log X}\frac{(\log_2 X)^{k-1}}{(k-1)!}.
\end{equation}
\end{theorem}

On peut aussi obtenir ce résultat lorsque $k\in\left[\delta\log_2 X, (\e-\delta)\log_2 X\right]$ pour $\delta>0$ fixé, ce que le Théorème~$\ref{petits intervalles}$ ne nous permet pas de faire, mais la borne inférieure sur $h$ est alors bien moins bonne. On obtient en effet la condition $h\geq\exp\left((\log X)^{r\log r-r+1+\epsilon}\right)$ où l'on a posé $r:=\frac{k}{\log_2 X}$ et $\epsilon>0$ est fixé. La relation asymptotique $\pi_{\left\lfloor\log_2 x\right\rfloor}(x)\sim \frac{x}{\sqrt{2\pi\log_2 x}}$ nous suggère que la condition $(\ref{borneinfh})$ pourrait être remplacée par $\frac{h}{\sqrt{\log_2 X}}\rightarrow +\infty$. Un travail en cours permettra de vérifier cela. Pour montrer les deux théorèmes, l'ingrédient principal est la proposition suivante.
\begin{prop}\label{propmatoradzi}
Uniformément pour tous $A,B>1$, pour toute fonction $\theta:\R\rightarrow\R$, pour tous $3\leq h\leq X$ et pour tout $\delta>0$, si l'on pose $I_{A,B}=\left[-A,-\frac{1}{B}\right]\cup\left[\frac{1}{B},A\right]$, alors
\begin{equation}\label{premierpoint}
\int_{I_{A,B}}\left\vert\frac{1}{h}\!\!\sum_{x<n\leq x+h}\!\!\e^{i\theta(\tau)\omega(n)\!\!}-\!\frac{1}{X}\!\!\sum_{X<n\leq 2X}\!\!\e^{i\theta(\tau)\omega(n)}\right\vert\frac{\mathrm{d}\tau}{\left\vert\tau\right\vert}\ll\left(\log (AB)\right)\left(\delta+\!\frac{\log_2 h}{\log h}\right)
\end{equation}
et\begin{equation}\label{secondpoint}
\int_{\left\vert z\right\vert=1}\left\vert\frac{1}{h}\!\!\sum_{x<n\leq x+h}\!\!z^{\omega(n)\!\!}-\!\frac{1}{X}\!\!\sum_{X<n\leq 2X}\!\!z^{\omega(n)}\right\vert\left\vert\mathrm{d}z\right\vert\ll\delta+\frac{\log_2h}{\log h}
\end{equation}
pour tout $x\in\left[X,2X\right]$ sauf sur un ensemble de mesure au plus
\[\ll X\left(\frac{\left(\log h\right)^{1/3}}{\delta^2h^{\delta/25}}+\frac{1}{\delta^2\left(\log X\right)^{1/50}}\right).\]
\end{prop}

Cette proposition, qui est montrée à la section~$\ref{sectionmatoradzi}$, découle d'un résultat récent et révolutionnaire de \MR\cite{matoradzi}. L'enjeu est de montrer que presque partout, la moyenne de la fonction arithmétique $n\mapsto \e^{i\theta\omega(n)}$ sur $(x,x+h]$ est proche de sa moyenne sur $[X,2X]$. Dans~\cite{matoradzi}, \MR traitent le cas des fonctions multiplicatives réelles, et précisent que le champ d'application peut-être étendu de diverses manières, notamment au cas des fonctions qui ne sont pas \og$p^{it}$-simulatrices\fg\ (\emph{pretentious} en anglais, notion développée depuis quelques années par Granville, Soundararajan, et d'autres auteurs). Ainsi, dans~\cite{matoradzitao}, Matomäki, Radziwi\l\l\ et Tao ont montré que la moyenne d'une fonction multiplicative complexe de module inférieur à $1$ et non \og$p^{it}$-simulatrice\fg\ était nulle sur presque tous les petits intervalles. Cependant, dans notre cas, lorsque $\theta$ tend vers $0$ la moyenne n'est plus nulle. On est donc amenés à comparer la moyenne sur les petits intervalles à celle sur un intervalle dyadique $\left[X,2X\right]$. Par ailleurs, les intégrales nous amènent à passer par une estimation plus forte sur un sous-ensemble $\mathcal{S}$ de $\left[X,2X\right]$, comme dans~\cite{matoradzi}.

On termine cette introduction en mentionnant un résultat récent de Teräväinen~\cite{teravainen}, qui a montré que pour tout $k\geq 3$ fixé et $\varepsilon>0$, presque tout intervalle $\left[x,x+\left(\log x\right)^{1+\varepsilon}\right]$ contient un entier $n$ tel que $\omega(n)=k$.

\section{Remerciements}

\indent\indent Je tiens à remercier chaleureusement mon directeur de thèse, Régis de la Bretèche, pour les nombreuses discussions que nous avons eues, et ses suggestions fort utiles. Je remercie par ailleurs Maksym Radziwi\l\l\ pour ses remarques pertinentes.

\section{Notations}

\indent\indent La lettre $p$ est réservée aux nombres premiers. Par conséquent, on écrit par exemple $\sum_p$ au lieu de $\sum_{p\text{ premier}}$. La lettre $\omega$ désigne la fonction additive qui compte le nombre de facteurs premiers distincts : $\omega(n)=\sum_{p\vert n}1$. Si $E$ est un ensemble fini, on note $\#E$ son cardinal. On écrit $f\ll g$ ou $f=\O{g}$ (resp. $f\gg g$) pour dire qu'il existe une constante absolue $C>0$ telle que $\left\vert f\right\vert\leq C\left\vert g\right\vert$ (resp. $\left\vert f\right\vert\geq C\left\vert g\right\vert$). La région de validité de cette inégalité, si elle n'est pas précisée, est claire d'après le contexte. On écrit par exemple $\ll_{\varepsilon}$ ou $O_{\varepsilon}$ pour signifier que la constante implicite $C$ peut dépendre de $\varepsilon$. La relation $f\asymp g$ signifie que l'on a simultanément $f\ll g$ et $f\gg g$. On utilise la notation (asymétrique) $a\sim b$ pour dire $b<a\leq2b$.

\section{Quelques lemmes classiques}

\indent\indent Pour montrer le Théorème~$\ref{petits intervalles}$, on reprend essentiellement la preuve du théorème d'Erd\o s-Kac qui figure dans~\cite[théorème~\MakeUppercase{\romannumeral 3}.4.15]{Tenenbaum}. Celle-ci repose principalement sur l'estimation de la somme des $\e^{i\theta\omega(n)}$ et sur l'inégalité de Berry-Esseen, qui permet de relier les fonctions de répartition à leur fonction caractéristique. Cependant, $\Phi_X$ n'est pas exactement une fonction de répartition, elle contient même des discontinuités (défaut qui a été introduit afin d'avoir un meilleur terme d'erreur), on utilise alors des travaux de la thèse d'Esseen, qui permettent de contourner le problème.

\begin{lemme}\label{BE}
Soit $m>0$ un réel fixé. Soit $F$ une fonction de répartition et $f$ sa fonction caractéristique définie par
\[f(\tau):=\int_{-\infty}^{+\infty}\e^{i\tau x}\mathrm{d}F(x).\]
Soit $G$ une fonction réelle à variation bornée sur $\R$, de fonction caractéristique $g$. On suppose les faits suivants vérifiés :
\begin{itemize}
\item $F(-\infty)=G(-\infty)$ et $F(+\infty)=G(+\infty)$,
\item si $G$ est discontinue en deux points $x$ et $y$ distincts, alors $\left\vert x-y\right\vert>m$,
\item G est dérivable en tout point de continuité, de dérivée bornée en valeur absolue par $\left\Vert G'\right\Vert_{\infty}$,
\item F ne peut être discontinue en $x$ que si $G$ l'est aussi.
\end{itemize}
Alors pour tout $T\gg\frac{1}{m}$, on a
\begin{equation}\label{inegBE}
\left\Vert F-G\right\Vert_{\infty}\ll\frac{\left\Vert G'\right\Vert_{\infty}}{T}+\int_{-T}^T\left\vert\frac{f(\tau)-g(\tau)}{\tau}\right\vert\mathrm{d}\tau.
\end{equation}
\end{lemme}
\begin{proof}
Voir~\cite[theorem~\MakeUppercase{\romannumeral 2}.2.b]{theseEsseen}.
\end{proof}
Avant de procéder à la démonstration du Théorème~$\ref{petits intervalles}$, on a besoin des deux lemmes suivants.
\begin{lemme}\label{TK}
Pour $X\geq 20$, on a
\[\frac{1}{X}\sum_{n\sim X}\left\vert\omega(n)-\log_2 X\right\vert\ll\sqrt{\log_2 X}.\]
\end{lemme}
\begin{proof}
Cette inégalité, qui est initialement due à Tur\'an, est une application directe de l'inégalité de Tur\'an-Kubilius, dont la démonstration se trouve par exemple dans~\cite[théorème~\MakeUppercase{\romannumeral 3}.3.1]{Tenenbaum}.
\end{proof}
\begin{lemme}\label{SD}
Uniformément pour $X\geq 3$ et $t\in\R$,
\[\frac{1}{X}\sum_{X<n\leq 2X}\e^{it\omega(n)}=A\left(\e^{it}\right)\left(\log X\right)^{\e^{it}-1}+\O{(\log X)^{\cos{t}-2}},\]
où $A(z)=1+c_1\left(z-1\right)+\O{(z-1)^2}$ dans la région $z-1=o(1)$, où $c_1$ est définie par $(\ref{mertensconstant})$.
\end{lemme}
\begin{proof}
Se déduit facilement de~\cite[théorème~\MakeUppercase{\romannumeral 2}.6.1]{Tenenbaum}. La constante $c_1$ apparaît comme la dérivée en $1$ de la fonction $z\mapsto\Gamma(z)^{-1}\prod_p\left(1+\frac{z}{p-1}\right)\left(1-\frac{1}{p}\right)^z$. Pour la calculer, on utilise la relation $\Gamma^\prime(1)=-\gamma$ (\emph{cf.}~\cite[corollaire~\MakeUppercase{\romannumeral 2}.0.7]{Tenenbaum}).
\end{proof}

\section{Démonstration du Théorème~\ref{petits intervalles}}

\indent\indent Soient $20\leq h\leq X$ et $x\sim X$. Majorons $\left\Vert F_{(x,x+h]}-\Phi_X\right\Vert_{\infty}$. On se restreint à $h\in\N$ pour des raisons qui apparaissent plus tard. Cela se fait au prix d'un $\O{\frac{1}{h}}$ qui est négligeable. Dans le but d'utiliser l'inégalité d'Esseen, on introduit les fonctions caractéristiques, lorsque $x\sim X$ et $\tau\in\R$,
\begin{align}
f_{(x,x+h]}(\tau) &:=\int_{-\infty}^{+\infty}\e^{i\tau y}\mathrm{d}F_{(x,x+h]}(y)=\frac{1}{h}\sum_{x<n\leq x+h}\e^{i\tau\frac{\omega(n)-\log_2 X}{\sqrt{\log_2 X}}},\\
\phi_X(\tau) &:=\int_{-\infty}^{+\infty}\e^{i\tau y}\mathrm{d}\Phi_X(y)=\e^{-\frac{\tau^2}{2}}\left(1+\frac{i\tau c_1-i\frac{\tau^3}{6}}{\sqrt{\log_2 X}}\right)+\Delta_X(\tau),\label{fourierPhi}
\end{align}
où
\begin{equation}\label{Delta}
\Delta_X(\tau):=\frac{-i\tau}{2\pi\sqrt{\log_2 X}}\sum_{\nu\in\Z^*}\frac{\e^{2i\pi\nu\log_2 X}}{i\nu}\e^{-\frac{1}{2}\left(\tau+2\pi\nu\sqrt{\log_2 X}\right)^2}.
\end{equation}
Le calcul pour $\phi_X(\tau)$ est aisé une fois connu le développement en série de Fourier
\[\frac{1}{2}-\langle t\rangle=\sum_{\nu\in\Z^*}\frac{\e^{2i\pi\nu t}}{2i\pi\nu}\hspace{2cm}\left(t\notin\N\right).\]
Pour simplifier les notations, on note désormais $T:=\log_2 X$. Soit $1<A\leq T$. Avec le Lemme $\ref{BE}$, on a alors
\begin{equation}\label{avecBE}
\left\Vert F_{(x,x+h]}-\Phi_X\right\Vert_{\infty}\ll\frac{1}{A}+\int_{-A}^A\left\vert\frac{f_{(x,x+h]}(\tau)-\phi_X(\tau)}{\tau}\right\vert\mathrm{d}\tau.
\end{equation}
L'intégrale est traitée en deux parties. On se donne $B>1$ et on étudie d'abord la contribution de l'intervalle $\left[-\frac{1}{B},\frac{1}{B}\right]$. Pour tout $y\in\R$, on a $\e^{iy}=1+\O{y}$. Alors pour tous $x\sim X$ et $\tau\in\R$, en se rappelant que $h$ est entier, on obtient
\[f_{(x,x+h]}(\tau)=1+\O{\frac{\left\vert\tau\right\vert}{\sqrt{T}h}\sum_{x<n\leq x+h}\left\vert\omega(n)-T\right\vert}.\]
Par ailleurs, pour $\left\vert\tau\right\vert\leq 1$ on a $\phi_X(\tau)=1+\O{\left\vert\tau\right\vert}$. Alors 
\[I(x):=\int_{-1/B}^{1/B}\left\vert\frac{f_{(x,x+h]}(\tau)-\phi_X(\tau)}{\tau}\right\vert\mathrm{d}\tau\ll\frac{1}{B}+\frac{1}{B\sqrt{T}h}\sum_{x<n\leq x+h}\left\vert\omega(n)-T\right\vert,\]
et donc
\[\frac{1}{X}\int_X^{2X}I(x)\mathrm{d}x\ll\frac{1}{B}+\frac{1}{B\sqrt{T}hX}\int_X^{2X}\sum_{x<n\leq x+h}\left\vert\omega(n)-T\right\vert\mathrm{d}x.\]
En intégrant par rapport à $x$, chaque $\left\vert\omega(n)-T\right\vert$ pour $n\in\left(X,2X+h\right]$ apparaît moins de $h$ fois. Ainsi, avec le Lemme~$\ref{TK}$, comme $h\leq X$ on obtient
\[\frac{1}{X}\int_X^{2X}I(x)\mathrm{d}x\ll\frac{1}{B}.\]
Ainsi, pour tout $\delta_1>0$, on a la majoration
\begin{equation}\label{I0}
I(x)\ll\frac{\delta_1}{B}
\end{equation}
pour tout $x\sim X$ sauf sur un ensemble de mesure au plus
\begin{equation}\label{excI0}
\ll \frac{X}{\delta_1}.
\end{equation}

Majorons maintenant la contribution des $\tau\in I_{A,B}:=\left[-A,-\frac{1}{B}\right]\cup\left[\frac{1}{B},A\right]$, à l'aide de~$(\ref{premierpoint})$ de la Proposition~$\ref{propmatoradzi}$. Pour cela, on introduit pour $y\in\R$
\begin{equation}
F_X(y):=\frac{1}{X}\#\left\{n\sim X\,:\quad\omega(n)\leq\log_2 X+y\sqrt{\log_2 X}\right\},
\end{equation}
et sa fonction caractéristique, pour $\tau\in\R$,
\begin{equation}
f_X(\tau):=\int_{-\infty}^{+\infty}\e^{i\tau y}\mathrm{d} F_X(y)=\frac{1}{X}\sum_{n\sim X}\e^{i\tau\frac{\omega(n)-\log_2 X}{\sqrt{\log_2 X}}}.
\end{equation}
Alors pour tout $x\sim X$,
\begin{equation}\label{reste}
\int_{I_{A,B}}\left\vert\frac{f_{(x,x+h]}(\tau)-\phi_X(\tau)}{\tau}\right\vert\mathrm{d}\tau\ll I_1(x)+I_2
\end{equation} où
\begin{align*}
I_1(x) &:=\int_{I_{A,B}}\left\vert f_{(x,x+h]}(\tau)-f_X(\tau)\right\vert\frac{\mathrm{d}\tau}{\left\vert\tau\right\vert},\\
I_2 &:=\int_{I_{T,B}}\left\vert\frac{f_X(\tau)-\phi_X(\tau)}{\tau}\right\vert\mathrm{d}\tau.
\end{align*}
puisque $A\leq T$. Avec~$(\ref{premierpoint})$ de la Proposition~$\ref{propmatoradzi}$, pour tout $\delta_2>0$ on a
\begin{equation}\label{I1}
I_1(x)\ll\left(\log (AB)\right)\left(\delta_2+\frac{\log_2 h}{\log h}\right)
\end{equation}
pour tout $x\sim X$ sauf sur un ensemble de mesure au plus
\begin{equation}\label{excI1}
\ll X\left(\frac{\left(\log h\right)^{1/3}}{\delta_2^2h^{\delta_2/25}}+\frac{1}{\delta_2^2\left(\log X\right)^{1/50}}\right).
\end{equation}
Par ailleurs, uniformément pour $\left\vert\tau\right\vert\leq T^{1/6}$, avec $(\ref{fourierPhi})$ et le Lemme~$\ref{SD}$ on a
\[\left\vert f_X(\tau)-\phi_X(\tau)\right\vert\ll \e^{-\frac{\tau^2}{2}}\frac{\left\vert\tau\right\vert^2+\left\vert\tau\right\vert^6}{T}+\frac{1}{\log X}.\]
De plus, pour tout $\left\vert t\right\vert\leq\pi$, on a $\cos t-1\leq-\frac{2}{\pi^2}t^2$. Donc toujours avec le Lemme~$\ref{SD}$, pour tout $T^{1/6}\leq\left\vert\tau\right\vert\leq\pi\sqrt{T}$, on a 
\begin{equation}
\left\vert f_X(\tau)-\phi_X(\tau)\right\vert\ll \exp\left(-\frac{2}{\pi^2}\tau^2\right).
\end{equation}
On a alors
\begin{equation}\label{I2}
I_2\ll\frac{1}{T}+\frac{\log (TB)}{\log X}+\int_{\pi\sqrt{T}}^T\left\vert\frac{f_X(\tau)-\Delta_X(\tau)}{\tau}\right\vert\mathrm{d}\tau.
\end{equation}
Il reste à traiter la dernière intégrale. La fonction $\Delta_X$ est concentrée autour des multiples entiers de $2\pi\sqrt{T}$. C'est aussi le cas de la fonction $f_X$. On découpe alors cette intégrale en fonction de l'entier en question. On veut donc majorer, pour $1\leq k\leq \frac{\sqrt{T}}{2\pi}$,
\[J_k:=\int_{2\pi\left(k-\frac{1}{2}\right)\sqrt{T}}^{2\pi\left(k+\frac{1}{2}\right)\sqrt{T}}\left\vert\frac{f_X(\tau)-\Delta_X(\tau)}{\tau}\right\vert\mathrm{d}\tau=\int_{-\pi\sqrt{T}}^{\pi\sqrt{T}}\left\vert\frac{f_X(2k\pi\sqrt{T}+u)-\Delta_X(2k\pi\sqrt{T}+u)}{2k\pi\sqrt{T}+u}\right\vert\mathrm{d}u.\]
Concernant la somme $\Delta_X$, seul le terme pour $\nu=-k$ contribue de manière non négligeable à $J_k$. Ce terme vaut, en $\tau=2k\pi\sqrt{T}+u$,
\[\frac{2k\pi\sqrt{T}+u}{2\pi\sqrt{T}}\frac{\e^{-2ik\pi T}}{k}\e^{-\frac{u^2}{2}}=\e^{-2ik\pi T}\e^{-\frac{u^2}{2}}\left(1+\frac{u}{2k\pi\sqrt{T}}\right).\]
Par ailleurs, $\tau\mapsto \e^{i\tau\sqrt{T}}f_X(\tau)$ est $2\pi\sqrt{T}$-périodique. On a donc $f_X(2k\pi\sqrt{T}+u)=\e^{-2ik\pi T}f_X(u)$. Ainsi, on peut écrire
\[J_k\ll\frac{1}{k\sqrt{T}}\int_{-\pi\sqrt{T}}^{\pi\sqrt{T}}\left\vert f_X(u)-\e^{-\frac{u^2}{2}}\right\vert\mathrm{d}u+\frac{1}{kT}+\e^{-T},\]
où le dernier terme provient de la somme des termes $\nu\neq -k$. L'intégrale ci-dessus se majore de la même manière que précédemment, et on obtient alors $J_k\ll\frac{1}{kT}+\frac{1}{\log X}$. En sommant et en utilisant avec~($\ref{I2}$), il vient alors
\begin{equation}\label{I2final}
I_2\ll\frac{\log T}{T}+\frac{T\log (TB)}{\log X}.
\end{equation}
Finalement, avec $(\ref{avecBE})$, $(\ref{I0})$, $(\ref{excI0})$, $(\ref{reste})$, $(\ref{I1})$, $(\ref{excI1})$ et $(\ref{I2final})$, on a pour tous $\delta_1,\delta_2>0$
\begin{equation}
\left\Vert F_{(x,x+h]}-\Phi_X\right\Vert_{\infty}\ll\frac{1}{A}+\frac{\log T}{T}+\frac{\delta_1}{B}+\left(\log (AB)\right)\left(\delta_2+\frac{\log_2 h}{\log h}\right)+\frac{T\log (TB)}{\log X}
\end{equation}
pour tout $x\sim X$ sauf sur un ensemble de mesure au plus
\begin{equation}
\ll X\left(\frac{1}{\delta_1}+\frac{\left(\log h\right)^{1/3}}{\delta_2^2h^{\delta_2/25}}+\frac{1}{\delta_2^2\left(\log X\right)^{1/50}}\right).
\end{equation}
Pour finir, il suffit de choisir $A,B,\delta_1$ et $\delta_2$ convenablement. On se donne $\alpha\geq 1$. Dans le cas $\alpha\geq\left(\log X\right)^{1-\frac{2}{150}}$, le résultat est trivial puisque $\log h\leq\log X$ et $\left\Vert F_{(x,x+h]}-\Phi_X\right\Vert_{\infty}\ll 1$. Dans le cas $\alpha\leq\left(\log X\right)^{1-\frac{2}{150}}$ et $\frac{\log h}{\alpha\log_2 h}\leq\left(\log X\right)^{\frac{1}{150}}$, on choisit 
\[A=\min\left(\log h, T\right),\ B=\left(\log h\right)^{10\alpha+1},\ \delta_1=\left(\log h\right)^{10\alpha},\ \delta_2=100\alpha\frac{\log_2 h}{\log h}.\]
Enfin, si $\frac{\log h}{\alpha\log_2 h}\geq\left(\log X\right)^{\frac{1}{150}}$, on choisit
\[A=T,\ B=\left(\log X\right)^2,\ \delta_1=\log X,\ \delta_2=100\left(\log X\right)^{-\frac{1}{150}}.\]
\qed

\section{Démonstration du Théorème~\ref{loislocales}}

\indent\indent Soit $k\geq 1$. Soient $20\leq h\leq X$ et $x\sim X$. Avec la formule de Cauchy, on a
\[\left\vert\pi_k(x+h)-\pi_k(x)-\frac{h}{X}\left(\pi_k(2X)-\pi_k(X)\right)\right\vert=\left\vert\frac{h}{2i\pi}\int_{\left\vert z\right\vert =1}\left(\frac{1}{h}\!\!\sum_{x<n\leq x+h}\!\!z^{\omega(n)\!\!}-\!\frac{1}{X}\!\!\sum_{X<n\leq 2X}\!\!z^{\omega(n)}\right)\frac{\mathrm{d}z}{z^{k+1}}\right\vert.\]
Ainsi, avec la majoration~$(\ref{secondpoint})$ de la Proposition~$\ref{propmatoradzi}$ lorsque $\delta=100\frac{\log_2h}{\log h}$, on obtient
\begin{equation}
\left\vert\pi_k(x+h)-\pi_k(x)-\frac{h}{X}\left(\pi_k(2X)-\pi_k(X)\right)\right\vert\ll h\frac{\log_2 h}{\log h}
\end{equation}
pour presque tout $x\sim X$ lorsque $h$ tend vers l'infini. Par ailleurs, avec~\cite[théorème~\MakeUppercase{\romannumeral 2}.6.4]{Tenenbaum}, uniformément pour $k=\log_2 X+\O{\sqrt{\log_2 X}}$, on a
\begin{equation}
\frac{\pi_k(2X)-\pi_k(X)}{X}\sim\frac{1}{\log X}\frac{(\log_2 X)^{k-1}}{(k-1)!}\asymp\frac{1}{\sqrt{\log_2 X}}.
\end{equation}
On en déduit le résultat recherché lorsque $h\geq\exp\left((\log_2 X)^{1/2+\epsilon}\right)$ pour un certain $\epsilon>0$ fixé.
\qed\\
Comme précisé dans l'introduction, la preuve fonctionne lorsque $k\in\left[\delta\log_2X, (\e-\delta)\log_2X\right]$ avec $\delta>0$ fixé, mais on obtient alors une moins bonne borne inférieure pour $h$.

\section{Estimation du type Matom\"{a}ki-Radziwi\l\l}\label{sectionmatoradzi}

\indent\indent La Proposition~$\ref{propmatoradzi}$ se déduit de la Proposition~$\ref{propsurS}$ ci-dessous, qui est une conséquence directe de~\cite{matoradzi}. On commence par énoncer la Proposition~$\ref{propsurS}$, dont la démonstration est donnée à la section~$\ref{adaptation}$, et on montre comment en déduire la Proposition~$\ref{propmatoradzi}$.

\subsection{Première simplification : définition de l'ensemble $\mathcal{S}$ et rétrécissement de l'intervalle de sommation}\label{sectionmatoradzi}

\indent\indent Soit $X\geq 20$. La Proposition~$\ref{propmatoradzi}$ découle d'une meilleure estimation sur un sous-ensemble dense $\mathcal{S}_X\subset[X,2X]$ défini de la manière suivante. Soit $\eta\in\left(0,\frac{1}{6}\right)$. On considère la suite d'intervalles $\left[P_j,Q_j\right]$ définis comme ceci :
\begin{itemize}
\item[$\bullet$] $1\leq\left(\log Q_1\right)^{40/\eta}\leq P_1\leq Q_1\leq \exp\left(\sqrt{\log X}\right)$,
\item[$\bullet$] pour $1\leq j\leq J$, $P_j=\exp\left(j^{4j}\left(\log Q_1\right)^{j-1}\left(\log P_1\right)\right)$,
\item[$\bullet$] pour $1\leq j\leq J$, $Q_j=\exp\left(j^{4j+2}\left(\log Q_1\right)^{j}\right)$,
\end{itemize}
où $J$ est défini comme le maximum des indices $j$ pour lequel $Q_j\leq\exp\left(\sqrt{\log X}\right)$. On a donc $J\ll\frac{\log_2 X}{\log_3 X}$. On définit alors $\mathcal{S}_X$ comme l'ensemble des entiers de $[X,2X]$ qui contiennent au moins un facteur premier dans chacun des $[P_j,Q_j]$ pour $1\leq j\leq J$. On remarque que
\[\frac{\log P_j}{\log Q_j}=\frac{1}{j^2}\frac{\log P_1}{\log Q_1},\hspace{1cm}(1\leq j\leq J)\]
et donc, par le lemme fondamental du crible, l'ensemble $[X,2X]\smallsetminus\mathcal{S}_X$ a une densité $\ll\frac{\log P_1}{\log Q_1}$.

\begin{prop}\label{propsurS}
Soit $\eta\in\left(0,\frac{1}{6}\right)$ fixé. Alors uniformément pour $2\leq h\leq h_2=\frac{X}{\left(\log X\right)^{1/5}}$ et $\theta\in\R$, si $\mathcal{S}=\mathcal{S}_X$ est défini comme ci-dessus avec $\left[P_1,Q_1\right]\subset\left[1,h\right]$, on a
\[\frac{1}{X}\int_X^{2X}\left\vert\frac{1}{h}\sum_{\substack{x<n\leq x+h \\ n\in\mathcal{S}}}\e^{i\theta\omega(n)}-\frac{1}{h_2}\sum_{\substack{x<n\leq x+h_2 \\ n\in\mathcal{S}}}\e^{i\theta\omega(n)}\right\vert^2\mathrm{d}x\ll\frac{\left(\log h\right)^{1/3}}{P_1^{1/6-\eta}}+\frac{1}{\left(\log X\right)^{1/50}}.\]
\end{prop}
Cette proposition est l'analogue de ce que l'on obtient en combinant~\cite[lemma~14]{matoradzi} et~\cite[proposition~1]{matoradzi}. Elle est démontrée à la section~$\ref{adaptation}$. On énonce maintenant deux lemmes qui nous permettent de vérifier que la Proposition~$\ref{propsurS}$ implique bien la Proposition~$\ref{propmatoradzi}$.
\begin{lemme}\label{raccourcirunpeu}
Soit $X\geq 3$. Uniformément pour $\theta\in\R$, $x\sim X$ et $\frac{X}{\left(\log X\right)^{1/5}}\leq y\leq X$, on a
\[\frac{1}{y}\sum_{x<n\leq x+y}\e^{i\theta\omega(n)}=\frac{1}{X}\sum_{n\sim X}\e^{i\theta\omega(n)}+\O{\left(\log X\right)^{\cos\theta-9/5}}.\]
\end{lemme}
\begin{proof}
Se déduit du Lemme~$\ref{SD}$.
\end{proof}
\begin{lemme}\label{crible}
Il existe une constante $c>0$ telle que pour tous $3\leq h_2\leq X$, $x\sim X$ et $3\leq P\leq Q\leq\exp\left(\sqrt{\log X}\right)$ on ait
\[\#\bigg\{x<n\leq x+h_2\,:\quad\bigg(n,\!\!\prod_{P< p\leq Q}\!\!p\bigg)\!=1\bigg\}=h_2\!\!\prod_{P< p\leq Q}\!\!\left(1-\frac{1}{p}\right)+\O{X\exp\left(-c\sqrt{\log X}\right)}.\]
\end{lemme}
\begin{proof}
C'est une conséquence directe de~\cite[theorem~6.1]{FIoperadecribro} avec par exemple $D=X^{1/2}$ et $z=\exp\left(\sqrt{\log X}\right)$.
\end{proof}
On peut désormais montrer la Proposition~$\ref{propmatoradzi}$.
\begin{proof}[Démontration de la Proposition~\ref{propmatoradzi}]
On ne détaille ici que l'obtention du premier point $(\ref{premierpoint})$, le second se traitant de manière tout à fait analogue. Soient $A,B>1$, $\delta>0$ et $\theta:\R\rightarrow\R$. Soient $3\leq h\leq X$. Si $h\geq h_2:=\frac{X}{\left(\log X\right)^{1/5}}$, le résultat est directement obtenu par l'application du Lemme~$\ref{raccourcirunpeu}$ puisque l'on peut imposer $\delta\gg\left(\log X\right)^{-1/99}$. On suppose donc désormais $h\leq h_2$.
Pour tout $x\sim X$, avec le Lemme~$\ref{raccourcirunpeu}$ on a
\begin{equation}
\int_{I_{A,B}}\left\vert\frac{1}{h_2}\sum_{x<n\leq x+h_2}\e^{i\theta(\tau)\omega(n)}-\frac{1}{X}\sum_{n\sim X}\e^{i\theta(\tau)\omega(n)}\right\vert\frac{\mathrm{d}\tau}{\left\vert\tau\right\vert}\ll\frac{\log (AB)}{\left(\log X\right)^{4/5}}.
\end{equation}
On prend $\mathcal{S}$ comme dans la Proposition~$\ref{propsurS}$ avec $\left[P_1,Q_1\right]\subset\left[1,h\right]$. On effectue la manipulation suivante (qui se trouve dans~\cite{matoradzi}) :
\begin{align*}
\frac{1}{h}\sum_{\substack{x<n\leq x+h \\ n\notin\mathcal{S}}}1 &=1-\frac{1}{h}\sum_{\substack{x<n\leq x+h \\ n\in\mathcal{S}}}1+\O{\frac{1}{h}}\\
&=\frac{1}{h_2}\sum_{\substack{x<n\leq x+h_2\\ n\in\mathcal{S}}}1-\frac{1}{h}\sum_{\substack{x<n\leq x+h\\ n\in\mathcal{S}}}1+\frac{1}{h_2}\sum_{\substack{x<n\leq x+h_2\\ n\notin\mathcal{S}}}1+O\left(\frac{1}{h}\right).
\end{align*}
Cela nous permet d'écrire, uniformément pour $\theta\in\R$,
\begin{align*}
&\left\vert\frac{1}{h}\sum_{\substack{x<n\leq x+h}}\e^{i\theta\omega(n)}-\frac{1}{h_2}\sum_{\substack{x<n\leq x+h_2}}\e^{i\theta\omega(n)}\right\vert \\
&\leq\left\vert\frac{1}{h}\sum_{\substack{x<n\leq x+h \\ n\in\mathcal{S}}}\e^{i\theta\omega(n)}-\frac{1}{h_2}\sum_{\substack{x<n\leq x+h_2 \\ n\in\mathcal{S}}}\e^{i\theta\omega(n)}\right\vert+\frac{1}{h}\sum_{\substack{x<n\leq x+h \\ n\notin\mathcal{S}}}1+\frac{1}{h_2}\sum_{\substack{x<n\leq x+h_2 \\ n\notin\mathcal{S}}}1\\
&\leq\left\vert\frac{1}{h}\!\!\sum_{\substack{x<n\leq x+h \\ n\in\mathcal{S}}}\!\!\!\!\e^{i\theta\omega(n)}\!-\frac{1}{h_2}\!\!\sum_{\substack{x<n\leq x+h_2 \\ n\in\mathcal{S}}}\!\!\!\!\e^{i\theta\omega(n)}\right\vert+\left\vert\frac{1}{h}\!\!\sum_{\substack{x<n\leq x+h \\ n\in\mathcal{S}}}\!\!\!\!\!1-\frac{1}{h_2}\!\!\sum_{\substack{x<n\leq x+h_2 \\ n\in\mathcal{S}}}\!\!\!\!\!1\,\right\vert+\frac{2}{h_2}\!\!\sum_{\substack{x<n\leq x+h_2\\ n\notin\mathcal{S}}}\!\!\!\!\!1+O\!\left(\frac{1}{h}\right)\!.
\end{align*}
D'après le Lemme~$\ref{crible}$, il existe une constante $c'>0$ telle que
\[\frac{1}{h_2}\sum_{\substack{x<n\leq x+h_2 \\ n\notin\mathcal{S}}}1\ll\sum_{j=1}^J\frac{\log P_j}{\log Q_j}+\exp\left(-c'\sqrt{\log X}\right)\ll\frac{\log P_1}{\log Q_1}.\]
On a alors
\begin{multline*}
\int_{I_{A,B}}\left\vert\frac{1}{h}\sum_{x<n\leq x+h}\e^{i\theta(\tau)\omega(n)}-\frac{1}{X}\sum_{n\sim X}\e^{i\theta(\tau)\omega(n)}\right\vert\frac{\mathrm{d}\tau}{\left\vert\tau\right\vert}\\
\ll\left(\log AB\right)\left(\frac{\log P_1}{\log Q_1}+\frac{1}{\left(\log X\right)^{4/5}}\right)+J_1(x)+J_2(x),
\end{multline*}
où
\begin{align*}
J_1(x):= &\int_{I_{A,B}}\left\vert\frac{1}{h}\sum_{\substack{x<n\leq x+h\\ n\in\mathcal{S}}}\e^{i\theta(\tau)\omega(n)}-\frac{1}{h_2}\sum_{\substack{x<n\leq x+h_2 \\ n\in\mathcal{S}}}\e^{i\theta(\tau)\omega(n)}\right\vert\frac{\mathrm{d}\tau}{\left\vert\tau\right\vert},\\
J_2(x):= &\int_{I_{A,B}}\left\vert\frac{1}{h}\sum_{\substack{x<n\leq x+h\\ n\in\mathcal{S}}}1-\frac{1}{h_2}\sum_{\substack{x<n\leq x+h_2 \\ n\in\mathcal{S}}}1\right\vert\frac{\mathrm{d}\tau}{\left\vert\tau\right\vert}.
\end{align*}
Par ailleurs,
\begin{align*}
&\frac{1}{X}\int_X^{2X}J_1(x)^2\mathrm{d}x=\int_{I_{A,B}}\int_{I_{A,B}}\frac{1}{\left\vert\tau\tau'\right\vert}\\
&\times\!\!\frac{1}{X}\!\!\int_X^{2X}\!\left\vert\frac{1}{h}\!\!\!\sum_{\substack{x<n\leq x+h\\ n\in\mathcal{S}}}\!\!\!\!\!\e^{i\theta(\tau)\omega(n)}\!-\!\frac{1}{h_2}\!\!\sum_{\substack{x<n\leq x+h_2 \\ n\in\mathcal{S}}}\!\!\!\!\!\e^{i\theta(\tau)\omega(n)}\right\vert\!\left\vert\frac{1}{h}\!\!\!\sum_{\substack{x<n\leq x+h\\ n\in\mathcal{S}}}\!\!\!\!\!\e^{i\theta(\tau')\omega(n)}\!-\!\frac{1}{h_2}\!\!\sum_{\substack{x<n\leq x+h_2 \\ n\in\mathcal{S}}}\!\!\!\!\!\e^{i\theta(\tau')\omega(n)}\right\vert\!\mathrm{d}x\\
&\times\mathrm{d}\tau\mathrm{d}\tau'.
\end{align*}
En appliquant successivement l'inégalité de Cauchy-Schwarz puis la Proposition~$\ref{propsurS}$ à l'intégrale par rapport à $x$, on obtient
\[\frac{1}{X}\int_X^{2X}J_1(x)^2\mathrm{d}x\ll\left(\log (AB)\right)^2\left(\frac{\left(\log h\right)^{1/3}}{P_1^{1/6-\eta}}+\frac{1}{\left(\log X\right)^{1/50}}\right).\]
Ainsi, pour tout $\delta>0$, on a
\[J_1(x)\leq\delta\log (AB)\]
pour tout $x\sim X$ sauf sur un ensemble de mesure au plus
\[\ll X\left(\frac{\left(\log h\right)^{1/3}}{\delta^2P_1^{1/6-\eta}}+\frac{1}{\delta^2\left(\log X\right)^{1/50}}\right).\]
L'intégrale $J_2(x)$ se traite exactement de la même façon (avec $\theta\left(\tau\right)=0$) et on trouve alors
\[\int_{I_{A,B}}\!\left\vert\frac{1}{h}\!\!\sum_{x<n\leq x+h}\!\!\e^{i\theta(\tau)\omega(n)}-\ \frac{1}{X}\!\sum_{n\sim X}\e^{i\theta(\tau)\omega(n)}\right\vert\frac{\mathrm{d}\tau}{\left\vert\tau\right\vert}\ll\left(\log (AB)\right)\left(\delta+\frac{\log P_1}{\log Q_1}+\!\frac{1}{\left(\log X\right)^{4/5}}\right)\]
pour tout $x\sim X$ sauf au plus
\[\ll X\left(\frac{\left(\log h\right)^{1/3}}{\delta^2P_1^{1/6-\eta}}+\frac{1}{\delta^2\left(\log X\right)^{1/50}}\right).\]
On peut omettre le terme $\left(\log X\right)^{-4/5}$ puisque l'on peut imposer $\delta\gg\left(\log X\right)^{-1/99}$.
Pour conclure, on fait le même choix pour $\eta, P_1$ et $Q_1$ que dans~\cite[section~9]{matoradzi} :
\begin{itemize}
\item[$\bullet$] $\displaystyle\eta=\frac{1}{150}$,
\item[$\bullet$] $\displaystyle Q_1=\min\left(h,\exp\left(\sqrt{\log X}\right)\right)$,
\item[$\bullet$] $\displaystyle P_1=\max\left(h^{\delta/4},\left(\log h\right)^{40/\eta}\right)$ si $\displaystyle h\leq\exp\left(\sqrt{\log X}\right)$, sinon $\displaystyle P_1=Q_1^{\delta/4}$.
\end{itemize}
\end{proof}

\subsection{Application du résultat de \MRf}\label{adaptation}
\indent\indent Il nous reste à montrer la Proposition~$\ref{propsurS}$. Comme indiqué précédemment, cette proposition est un application directe de~\cite[lemma~14]{matoradzi} et~\cite[proposition~1]{matoradzi} dans le cas où $f(n)=\e^{i\theta\omega(n)}$ ; cependant, la proposition 1 de~\cite{matoradzi} telle quelle, ne s'applique qu'aux fonctions multiplicatives $f:\N\rightarrow[-1,1]$ réelles. L'unique point de leur preuve ne couvrant pas le cas $f(n)=\e^{i\theta\omega(n)}$ se trouve page $28$ de~\cite{matoradzi}, c'est l'estimation 
\[\max_{\left(\log X\right)^{1/15}\leq\left\vert u\right\vert\leq 2T^{1+\varepsilon}}\left\vert R_{v,H}\left(1+iu\right)\right\vert\ll\left(\log X\right)^{-1/16+o(1)}\cdot\frac{\log Q}{\log P}\]
où
\begin{equation}\label{defRvH}
R_{v,H}(s):=\sum_{\substack{X\e^{-v/H}\leq n\leq 2X\e^{-v/H} \\ n\in\mathcal{S}}}\frac{f(n)}{n^s}\cdot\frac{1}{\#\left\{p\in\left[P,Q\right]\,:\quad p\vert n\right\}+1},\hspace{1cm}(s\in\C)
\end{equation}
qui est principalement due au fait qu'une fonction réelle n'est pas \og $p^{it}$-simulatrice\fg~\cite[lemma~2]{matoradzi}. Ce type de résultat est dû à des travaux initiés par Hal\'asz dans les années $1970$, puis amplifiés notamment par Granville et Soundararajan. Les fonctions qui nous intéressent ne sont pas \og $p^{it}$-simulatrices\fg\ non plus puisqu'elles sont constantes sur les nombres premiers. Cela nous permet d'adapter la preuve assez simplement. Pour la suite, on se donne l'ensemble $\mathcal{S}=\mathcal{S}_X$ défini précédemment et on pose $P:=\exp\left(\left(\log X\right)^{1-1/48}\right)$, $Q:=\exp\left(\frac{\log X}{\log_2 X}\right)$ et $H:=\left(\log X\right)^{1/48}$.
\begin{prop}\label{finale}
On suppose $\mathcal{S}, P, Q$ et $H$ définis comme ci-dessus. Soit $\varepsilon\in\left(0,1\right)$. Uniformément pour $3\leq T\leq X$, $\theta\in\R$ et $v\in\left[\left\lfloor H\log P\right\rfloor,H\log Q\right]$, lorsque $f(n)=\e^{i\theta\omega(n)}$ on a
\[\max_{\left(\log X\right)^{1/15}\leq\left\vert u\right\vert\leq 2T^{1+\varepsilon}}\left\vert R_{v,H}\left(1+iu\right)\right\vert\ll\left(\log X\right)^{-1/16+o(1)}\cdot\frac{\log Q}{\log P},\]
où $R_{v,H}$ est défini par $(\ref{defRvH})$, lorsque $X$ tend vers l'infini.
\end{prop}
La preuve suit le même schéma que dans~\cite[lemma~$3$]{matoradzi}. Lorsque $f$ et $g$ sont deux fonctions multiplicatives de module inférieur à $1$ et $x\geq 1$, on note
\[\mathbb{D}\left(f,g,x\right)^2=\sum_{p\leq x}\frac{1-\re\left(f(p)\overline{g(p)}\right)}{p}.\]
C'est au sens de cette \og distance\fg\ que l'on montre que nos fonctions ne ressemblent pas à $n\mapsto n^{it}$. En effet, on a le lemme suivant.
\begin{lemme}
Soient $f$ une fonction multiplicative à valeur dans le disque unité, et
\[F(s)=\sum_{n\sim x}\frac{f(n)}{n^s}.\hspace{1cm}(x\geq 1, s\in\C)\]
Uniformément pour $T_0\geq 1$, $t\in\R$ et $x\geq 1$, on a
\[\left\vert F\left(1+it\right)\right\vert\ll \left(1+m\right)\e^{-m}+\frac{1}{T_0}\]
où
\[m=m\left(x,T_0\right):=\min_{\left\vert t_0\right\vert\leq T_0}\mathbb{D}\left(f,n^{it+it_0},x\right)^2.\]
\end{lemme}
\begin{proof}
Se déduit directement de~\cite[corollaire~\MakeUppercase{\romannumeral 3}.4.12]{Tenenbaum} par sommation d'Abel.
\end{proof}
Le lemme suivant permet de contourner la condition non multiplicative $n\in\mathcal{S}$ --- on rappelle que l'ensemble $\mathcal{S}=\mathcal{S}_X$ a été défini au début de la section~\ref{sectionmatoradzi}.
\begin{lemme}\label{condSnonmult}
Pour $\mathcal{J}\subset\left\{1, ..., J\right\}$, on définit la fonction totalement multiplicative $g_{\mathcal{J}}$ par
\[g_{\mathcal{J}}(p)=\left\{
\begin{tabular}{ll}
1 & \text{si }$p\notin\bigcup_{j\in\mathcal{J}}\left[P_j,Q_j\right]$ \\
0 & \text{sinon.} 
\end{tabular}
\right.\]
Alors pour toute suite $\left(a_n\right)_{n\sim X}$, on a
\[\sum_{\substack{n\sim X \\ n\in\mathcal{S}}}a_n=\sum_{n\sim X}a_n\prod_{j=1}^J\left(1-g_{\left\{ j\right\}}(n)\right)=\sum_{\mathcal{J}\subset\left\{1, ..., J\right\}}\left(-1\right)^{\#\mathcal{J}}\sum_{n\sim X}a_ng_{\mathcal{J}}(n).\]
\end{lemme}
Avec ce lemme, on peut se passer de la condition $n\in\mathcal{S}$ au prix d'un facteur multiplicatif $2^J=\left(\log X\right)^{o(1)}$.
\begin{lemme}
Soient $A\geq 1$ et $\varepsilon\in\left(0,\frac{1}{4}\right)$ fixés. Soit $X\geq 20$. On suppose $\mathcal{S}=\mathcal{S}_X$, $P$ et $Q$ définis comme précédemment. Alors uniformément pour $\mathcal{J}\subset\left\{1, ..., J\right\}$, $\theta\in\R$, $1\leq\alpha\leq X^A$ et $x\in\left[X^{1/4},X\right]$, on a
\[\mathbb{D}\left(h_{\mathcal{J},\theta},n^{i\alpha},x\right)^2\geq\left(\frac{1}{3}-\frac{1}{48}-\varepsilon\right)\log_2 x+\O{1},\]
où $h_{\mathcal{J},\theta}(n)=g_{\mathcal{J}}(n)\e^{i\theta\omega(n)}$ lorsque $n$ n'a aucun facteur premier dans l'intervalle $\left[P,Q\right]$ et $h_{\mathcal{J},\theta}(n)=0$ sinon, avec les fonctions $g_{\mathcal{J}}$ définies au Lemme~$\ref{condSnonmult}$.
\end{lemme}
\begin{proof}
Puisque $Q_J\leq\exp\left(\sqrt{\log X}\right)$, on a
\begin{align*}
\mathbb{D}\left(h_{\mathcal{J},\theta},n^{i\alpha},x\right)^2 &=\sum_{p\leq x}\frac{1-\re\left(h_{\mathcal{J},\theta}(p)p^{-i\alpha}\right)}{p}\\
&\geq\sum_{\exp\left(\left(\log x\right)^{2/3+\varepsilon}\right)<p\leq\exp\left(\left(\log x\right)^{1-1/48}\right)}\frac{1-\re\left(\e^{i\theta}p^{-i\alpha}\right)}{p}\\
&\geq\left(\frac{1}{3}-\frac{1}{48}-\varepsilon\right)\log_2 x+\O{1}\\
& +\O{\sum_{\exp\left(\left(\log x\right)^{2/3+\varepsilon}\right)<p\leq\exp\left(\left(\log x\right)^{1-1/48}\right)}\frac{\re\left(\e^{i\theta}p^{-i\alpha}\right)}{p}}. &
\end{align*}
Par ailleurs, $\re\left(\e^{i\theta}p^{-i\alpha}\right)=\left(\cos\theta\right)\re\left(p^{-i\alpha}\right)-\left(\sin\theta\right)\im\left(p^{-i\alpha}\right)$. De plus, la région sans zéro de Korobov-Vinogradov pour la fonction zêta nous fournit (voir par exemple ligne $\left( 20\right)$ de~\cite{balazardtenenbaum})
\[\left\vert\sum_{\exp\left(\left(\log X\right)^{2/3+\varepsilon}\right)<p\leq\exp\left(\left(\log X\right)^{1-1/48}\right)}\frac{1}{p^{1+i\alpha}}\right\vert\ll 1.\]
On a bien le résultat attendu.
\end{proof}
\begin{proof}[Démonstration de la Proposition~\ref{finale}]
Grâce aux trois lemmes précédents, il suffit  de reprendre la preuve de~\cite[lemma~3]{matoradzi}, \emph{mutatis mutandis}.
\end{proof}

\bibliographystyle{plain-fr}
\bibliography{Bibliographie}

\textit{E-mail :} \url{eliegoudout@hotmail.com}
\end{document}